\documentclass[11pt,pdf]{article}
\usepackage{amsthm}
\usepackage{amssymb}
\usepackage{amsbsy}
\usepackage{amsfonts}
\usepackage{amsmath}
\usepackage{amsopn}

\newtheorem{theo}{Theorem}
\newtheorem{prop}[theo]{Proposition}
\newtheorem{lemm}[theo]{Lemma}

\DeclareMathOperator{\nr}{\mathbb{N}}

\DeclareMathOperator{\re}{\mathbb{R}}
\DeclareMathOperator{\p}{\mathbb{P}}
\DeclareMathOperator{\e}{\mathbb{E}}

\DeclareMathOperator{\vat}{\xrightarrow[t\to\infty]{}}

\begin{document}
\author{Jan Ob\l \'oj}
\date{}
\title{A complete characterization of local martingales which are functions of Brownian motion and its maximum}
\maketitle
\begin{center}\vspace{-0.8cm}
\textit{Laboratoire de Probabilit\'es et Mod\`eles Al\'eatoires, Universit\'e Paris 6}\\\textit{ 4 pl. Jussieu - Bo\^{i}te 188, 75252 Paris Cedex 05, France.}\\\textit{Department of Mathematics Warsaw University}\\ \textit{Banacha 2, 02-097 Warszawa, Poland.}\\\textit{e-mail: }\textrm{obloj@mimuw.edu.pl}
\end{center}
\begin{abstract}
We prove the max-martingale conjecture given in Ob\l \'oj and Yor \cite{ja_yor_maxmart}.
We show that for a continuous local martingale $(N_t:t\ge 0)$ and a function $H:\re\times\re_+\to \re$, $H(N_t,\sup_{s\leq t}N_s)$ is a local martingale if and only if there exists a locally integrable function $f$ such that
 $H(x,y)=\int_0^y f(s)ds-f(y)(x-y)+H(0,0)$. This implies readily, via L\'evy's equivalence theorem, an analogous result with the maximum process replaced by the local time at $0$.
\end{abstract}
\section{Introduction}
In our recent paper with Marc Yor \cite{ja_yor_maxmart}, we argued about the importance of a class of local martingales which are functions of the couple: continuous local martingale and its one-sided maximum process. We called them \emph{max-martingales} or simply \emph{M-martingales}. Such processes were first introduced by Az\'ema and Yor \cite{MR82c:60073a}, who described a family of such martingales, which is often referred to as \textit{Az\'ema-Yor martingales}. With Marc Yor \cite{ja_yor_maxmart}, we gave a complete description of this family. These martingales have a remarkably simple form, yet they proved to be a useful tool in various problems. We assembled with Marc Yor \cite{ja_yor_maxmart} applications including the Skorokhod embedding problem (cf.\ Ob\l \'oj \cite{genealogia}), a simple proof of Doob's maximal and $L^p$- inequalities, as well as a derivation of bounds on the possible laws of the maximum or the local time at $0$ of a continuous, uniformly integrable martingale, and links with Brownian penalization problems (cf.\ Roynette, Vallois and Yor \cite{rvy2}).

In this paper, we obtain a complete characterization of the max-martingales, which was conjectured in \cite{ja_yor_maxmart}. More precisely we show that Az\'ema-Yor martingales are actually the only continuous local martingales which are functions of Brownian motion and of its unilateral maximum process. Put differently \emph{Az\'ema-Yor martingales} are the only \emph{max-martingales} and we will use both terms interchangeably.
This paper is the first step in a more general project, which we hope to develop, and which consists in describing explicitly families of martingales which are functions of Brownian motion and its maximum, minimum and local time processes.

Section \ref{sec:main} below contains our main theorem, its corollaries and a complementary result. All proofs are gathered in the subsequent Section \ref{sec:proofs}. Section \ref{sec:optional} contains some arguments based on the optional stopping theorem which are very different from the arguments used in the proofs in Section \ref{sec:proofs}, and hopefully will give the reader some additional insight.

\section{Main results}
\label{sec:main}
Throughout, $N=(N_t:t\ge 0)$ denotes a continuous local martingale with $N_0=0$ and $\langle N\rangle_\infty=\infty$ a.s. Extensions of our results to local martingales with arbitrary $N_0$ are immediate. The maximum and minimum processes are denoted respectively: $\overline{N}_t=\sup_{s\le t} N_s$ and $\underline{N}_t=-\inf_{s\le t} N_s$, and $L^N_t$ denotes the local time at zero. Filtrations considered are always taken completed and right-continuous. 
$B=(B_t:t\ge 0)$ denotes a real-valued Brownian motion.
The following theorem is the main result of this paper.
\begin{theo}
\label{thm:maxmart} Let $\mathbf{D}=\{(x,y)\in \re^2: y\geq x\lor 0\}$ and $H:\mathbf{D}\to \re$ be a Borel function, and $N=(N_t:t\ge 0)$ be a continuous local martingale with $N_0=0$ and $\langle N\rangle_\infty=\infty$ a.s.
Then
$\big(H(N_t,\overline{N}_t):t\ge 0\big)$ is a right-continuous local martingale, in the natural filtration of $N$, if and only if
there exists $f:\re_+\to\re$, a locally integrable function, such that a.s., for all $t\ge 0$,
\begin{eqnarray}
  \label{eq:Hgenmart}
H(N_t,\overline{N}_t)&=&F(\overline{N}_t)-f(\overline{N}_t)(\overline{N}_t-N_t)+H(0,0)\\
&=&\int_0^t f(\overline{N}_s)dN_s+H(0,0), \label{eq:Hgenmartint}
\end{eqnarray}
where $F(y)=\int_0^y f(x)dx$.
\\Furthermore, if 
$f\geq 0$ and $\int_0^\infty f(x)dx<\infty$, then $\big(H(N_t,\overline{N}_t):t\ge 0\big)$ given in (\ref{eq:Hgenmart}) converges a.s., as $t\to \infty$, to $F(\infty)+H(0,0)$. If moreover $(N_t:t\ge 0)$ is a martingale with $\e \sup_{s\le t}|N_s|<\infty$, $t>0$, then the local martingale in (\ref{eq:Hgenmart}) is a martingale.
\end{theo}
\noindent Note that, in particular, if $H(N_t,\overline{N}_t)_{t\ge 0}$ is a right-continuous local martingale, then it is in fact a continuous local martingale and $H(\cdot,y)$ is a linear function for almost all $y>0$. We also can specify the maximum process of this local martingale. Indeed, from (\ref{eq:Hgenmart}) it is clear that $\sup_{s\le t}H(N_s,\overline{N}_s)=F(\overline{N}_t)+H(0,0)$.\\
The martingale property announced in the theorem was observed by Roynette, Vallois and Yor \cite{rvy2}.
We stress here that the local martingales we obtain have some interesting properties. If $f\ge 0$ and $F(\infty)<\infty$ then the process $M^f_t=F(\overline{N}_t)-f(\overline{N}_t)(\overline{N}_t-N_t)$ provides an example of a local martingale which converges a.s.\ to its maximum: $M^f_t\vat F(\infty)=\overline{M^f}_\infty$. Equivalently, $f(\overline{N}_t)(\overline{N}_t-N_t)$ is a local submartingale, zero at zero, which converges a.s.\ to zero as $t\to\infty$. 
\\ Theorem \ref{thm:maxmart} tell us that a local martingale $H(N_t,\overline{N}_t)$ is entirely characterized by its initial value $H(0,0)$ and by a locally integrable function $f$ such that (\ref{eq:Hgenmart}) holds. We point out that we can recover this function from the process $H(N_t,\overline{N}_t)$ itself. Indeed from (\ref{eq:Hgenmartint}) we have that $\langle N, H(N,\overline{N})\rangle_t=\int_0^t f(\overline{N}_s)d\langle N\rangle_s$. Thus the measure $d\langle N, H(N,\overline{N})\rangle_t$ is absolutely continuous with respect to $d\langle N\rangle_t$, and the density is given by $f(\overline{N}_t)$. This yields
\begin{equation}\label{eq:funrecov}
    f(x)=\frac{d\langle N, H(N,\overline{N})\rangle_{t}}{d\langle N\rangle_{t}}\bigg|_{t=T_x},
\end{equation}
where $T_x=\inf\{t:N_t=x\}$, $x>0$.\\
\noindent Naturally, Theorem \ref{thm:maxmart} has an analogue with the maximum $\overline{N}_t$ replaced by the minimum $\underline{N_t}$. More precisely, $H(N_t,\underline{N}_t)_{t\ge 0}$ is a right-continuous local martingale, if and only if, there exists a locally integrable function $f$ such that $H(N_t,\underline{N}_t)=F(\underline{N}_t)-f(\underline{N}_t)(\underline{N}_t+N_t)$ a.s. This is obtained upon substituting $N$ with $-N$ in Theorem \ref{thm:maxmart}.

Thanks to L\'evy's equivalence theorem, we can rephrase Theorem \ref{thm:maxmart} also in terms of the local time at zero instead of the maximum process.
\begin{theo}
  \label{thm:Lharm}
  Let $H:\re_+\times\re_+\to \re$ be a Borel function
and $N=(N_t:t\ge 0)$ be a continuous local
martingale with $N_0=0$ and $\langle N\rangle_\infty=\infty$ a.s. 
Then $\big(H(|N_t|,L^N_t):t\ge 0\big)$ is a right-continuous local martingale, in the natural filtration of $N$, if and only if there exists $g:\re_+\to\re$ a locally integrable function, such that a.s., for all $t\ge 0$,
\begin{eqnarray}
  \label{eq:HgenmartL}
H(|N_t|,L^N_t)&=&G(L^N_t)-g(L^N_t)|N_t|+H(0,0)\\
&=&-\int_0^t g(L^N_s)sgn(N_s)dN_s+H(0,0),  \label{eq:HgenmartLint}
\end{eqnarray}
 where $G(y)=\int_0^y g(x)dx$.
\\Furthermore, if 
$g\geq 0$ and $\int_0^\infty g(x)dx<\infty$, then $\big(H(|N_t|,L^N_t):t\ge 0\big)$ given in (\ref{eq:Hgenmart}) converges a.s., as $t\to \infty$, to $G(\infty)+H(0,0)$. If moreover $(N_t:t\ge 0)$ is a martingale with $\e \sup_{s\le t}|N_s|<\infty$, $t>0$, then the local martingale in (\ref{eq:Hgenmart}) is a martingale.
\end{theo}
\noindent Note that if $H(|N_t|,L^N_t)$ is a continuous local martingale then we can recover from it, in a similar manner to (\ref{eq:funrecov}), the function $g$ such that (\ref{eq:HgenmartL}) holds.

Theorem \ref{thm:maxmart} allows us also to consider local martingales of the form $H(N^+_t,\overline{N}_t)$, where $H$ is a Borel function. Indeed, $H(N^+_t,\overline{N}_t)$ can be written as $G(N_t,\overline{N}_t)$ with $G(x,y)=H(x\lor 0,y)$, and we can then apply Theorem \ref{thm:maxmart}.
This yields the following theorem.
\begin{theo}
\label{cor:martbplus}
Let $H:\re_+\times\re_+\to\re$ be a Borel function. Let $(N_t:t\ge 0)$ be a continuous local martingale with $N_0=0$ and $\langle N \rangle_\infty=\infty$ a.s.\ and put $N_t^+=\max\{N_t,0\}$, $N_t^-=\max\{-N_t,0\}$ and $N^*_t=\sup_{s\le t}|N_s|$.
Then the processes $(H(N_t^+,\overline{N}_t):t\ge 0)$, $(H(N_t^-,\underline{N}_t):t\ge 0)$ and $(H(|N_t|,N^*_t):t\ge 0)$ are right-continuous local martingales, in the natural filtration of $N$ if and only if they are a.s.\ constant.
\end{theo}
\noindent Note that it is not true that all martingales in the natural filtration of $N^+$ are also martingales in the natural filtration of $N$. In fact the former admits also discontinuous martingales unlike the latter.
From the proof it will be clear that the theorem stays true if we replace $N_t^+=\max\{N_t,0\}$ with some more complicated, appropriate function of $N_t$ and work with local martingales in the natural filtration of $N$.
For example we can easily see that if
 $A$ is an interval and $(-\infty,0)\nsubseteq A$, then $(H(N_t\mathbf{1}_{N_t\in A},\overline{N}_t):t\ge 0)$, is a right-continuous local martingale, in the natural filtration of $N$, if and only if it is constant a.s.\\ Similar reasonings can be naturally developed for the local time $L^N$ in the place of the maximum $\overline{N}$.

Finally we present a deterministic description of $H$'s such that (\ref{eq:Hgenmart}) holds. This is very close to studying the fine topology for the process $(B_t,\overline{B}_t)$.
\begin{prop}
\label{prop:Hgen}
In the setup of Theorem \ref{thm:maxmart}, (\ref{eq:Hgenmart}) holds
if and only if there exists a set $\Gamma\subset \mathbf{D}$ of Lebesgue measure zero, such that
  \begin{equation}
    \label{eq:Hgen}
 H(x,y)=F(y)-f(y)(y-x)+H(0,0),\quad\forall\, (x,y)\in \mathbf{D}\setminus \Gamma,
  \end{equation}
and $\Gamma_2=\{y:\exists x,\,(x,y)\in\Gamma\}$ has Lebesgue measure zero, and $\{(y,y):y\ge 0\}\cap \Gamma=\emptyset$.
\end{prop}
An analogue for the function $G$ satisfying (\ref{eq:HgenmartL}) follows.
\section{Proofs}
\label{sec:proofs}
We note that it suffices to prove Theorems \ref{thm:maxmart}, \ref{thm:Lharm}, \ref{cor:martbplus} and Proposition \ref{prop:Hgen} for $N=B$, a standard real-valued Brownian motion, as
then through Dambis-Dubins-Schwarz theorem (cf.\ Revuz and Yor \cite{MR2000h:60050} p.\ 181), it extends to any
continuous local martingale $N$ with $N_0=0$ and $\langle N\rangle_\infty=\infty$.\\ Indeed we know that if $T_u$ is the right-continuous inverse of $\langle N\rangle _t$ then the process $\beta_u=N_{T_u}$ is a Brownian motion and $N_t=\beta_{\langle N\rangle_t}$. It follows that $\overline{N}_t=\overline{\beta}_{\langle N\rangle_t}$.\\
Note $(\mathcal{F}^N_t)$ the natural filtration of $N$ and $\mathcal{G}_u=\mathcal{F}^N_{T_u}$, $\mathcal{F}^\beta_u=\sigma(\beta_s:s\le u)$ two filtrations with respect to which $(\beta_u)$ is a Brownian motion. Naturally $\mathcal{F}^\beta_u\subset\mathcal{G}_u$ 
but in fact the smaller filtration is immersed in the larger, meaning that all $(\mathcal{F}^\beta_u)$-local martingales are also $(\mathcal{G}_u)$-local martingales.
This follows readily from the representation of $(\mathcal{F}^\beta_u)$-local martingales as stochastic integrals with respect to $\beta$ and thus $(\mathcal{G}_u)$-local martingales (cf.\ Yor \cite{MR528910}).
This entails that $H(\beta_u,\overline{\beta}_u)$, which is $(\mathcal{F}^\beta_u)$-measurable, is a $(\mathcal{G}_u)$-local martingale if and only if it is also a $(\mathcal{F}^\beta_u)$-local martingale.\\
Thus if $H(\beta_u,\overline{\beta}_u)$ is a $(\mathcal{F}^\beta_u)$-local martingale, it is a $(\mathcal{G}_u)$-local martingale and therefore its time-changed version $H(\beta_{\langle N\rangle_t},\overline{\beta}_{\langle N\rangle_t})=H(N_t,\overline{N}_t)$ is a $(\mathcal{F}^N_t)$-local martingale (note that the time change is continuous).
Conversely, as $N$ is constant on the jumps of $(T_u)$, we have $\overline{\beta}_u=\overline{N}_{T_u}$ and if $H(N_t,\overline{N}_t)$ is a $(\mathcal{F}^N_t)$-local martingale, then its time-changed version $H(\beta_u,\overline{\beta}_u)$ is a $(\mathcal{G}_u)$-local martingale (cf.\ Revuz and Yor \cite{MR2000h:60050} Proposition V.1.5) and thus a $(\mathcal{F}^\beta_u)$-local martingale.\\
Likewise, as $L_t^N=\lim_{\epsilon\to 0}\frac{1}{2\epsilon}\int_0^t \mathbf{1}_{|N_s|\leq \epsilon}d\langle N\rangle_s$ (cf.\ Revuz and Yor \cite{MR2000h:60050} p.\ 227), it is easy to see that $L^B_{\langle M\rangle_t}=L^N_t$, $L^B_t=L^N_{T_t}$ and that Theorem \ref{thm:Lharm} holds for an arbitrary continuous local martingale $N$, with $N_0=0$ and $\langle N\rangle_\infty=\infty$ a.s., if and only if it holds for Brownian motion.
\\
Theorem \ref{thm:Lharm} for Brownian motion, follows from Theorem \ref{thm:maxmart} with
L\'evy's equivalence theorem, which grants that the processes
$((B_t,\overline{B}_t):t\geq 0)$ and $((L^B_t-|B_t|,L^B_t):t\geq 0)$
have the same distribution.

We turn to the proof of Theorem \ref{cor:martbplus}. Suppose that $(H(B^+_t,\overline{B}_t):t\ge 0)$ is a right-continuous local martingale in the natural filtration of $B$. Then $(G(B_t,\overline{B}_t):t\ge 0)$ is also a right-continuous local martingale, where $G(x,y)=H(x\lor 0,y)$. By Theorem \ref{thm:maxmart}, there exists a locally integrable function $f$ such that $G(x,y)=F(y)-f(y)(y-x)+G(0,0)$ a.e., where $F(y)=\int_0^y f(x)dx$. However, for any fixed $y\in \re_+$, $G$ is constant for $x<0$ which means that $f(x)=0$ a.e.\ (cf.\ proof of Proposition \ref{prop:Hgen} below) and thus $H(N^+_t,\overline{N}_t)=H(0,0)$ a.s. An analogous result for $H(N^-_t,\underline{N}_t)$ follows.\\
Let $A^+_t=\int_0^t\mathbf{1}_{B_s\ge 0}ds$ and $\alpha^+_u$ be its right-continuous inverse. Then $W_u=B^+_{\alpha^+_u}$ is a reflected Brownian motion in the filtration $\mathcal{G}_u=\mathcal{F}_{\alpha_u^+}$, where $(\mathcal{F}_t)$ is the natural filtration of $B$, and $\overline{W}_u=\overline{B}_{\alpha^+_u}$. If we write $(\mathcal{F}^W_u)$ the natural filtration of $W$, then $H(W_u,\overline{W}_u)$ is a $(\mathcal{G}_u)$-local martingale if and only if it is also a $(\mathcal{F}^W_u)$-local martingale. This follows from our discussion above and the fact that $W$ can be written as $W_u=\beta_u+L_u^\beta$, where $\beta_u$ is a $(\mathcal{G}_u)$ Brownian motion and the natural filtrations of $W$ and $\beta$ are equal (cf. Yor \cite{MR0471060}).\\
Suppose now that $H(W_u,\overline{W}_s)$ is a $(\mathcal{F}_t^W)$ right-continuous local martingale and thus a $(\mathcal{G}_t)$ right-continuous local martingale. As the time change $A^+_t$ is continuous, the time-changed version $H(W_{A_t^+},\overline{W}_{A^+_t})=H(B^+_t,\overline{B}_t)$ is a $(\mathcal{F}_t)$ right-continuous local martingale and thus is constant a.s. This completes the proof of Theorem \ref{cor:martbplus}.

We now prove Proposition \ref{prop:Hgen}. As the law of $(B_t,\overline{B}_t)$ is equivalent to the Lebesgue measure on $\mathbf{D}$, it is clear that (\ref{eq:Hgenmart}) implies that (\ref{eq:Hgen}) holds for $(x,y)\in \mathbf{D}\setminus \Gamma$, for some set $\Gamma$ of Lebesgue measure zero. However we know also that $H(B_t,\overline{B}_t)$ is a.s.\ continuous and this yields more constraints on $\Gamma$. Conversely, if (\ref{eq:Hgen}) holds for a Lebesgue null set $\Gamma$, small enough so that $H(B_t,\overline{B}_t)$ is a.s. continuous, then (\ref{eq:Hgenmart}) holds.\\
For a set $\Gamma\subset \mathbf{D}$ let $\Gamma_2:=\{y:\exists x,\,(x,y)\in \Gamma\}$, and for $y\in \Gamma_2$, $x_y:=\sup\{x:(x,y)\in \Gamma\}$.
Note $\Gamma_2^+=\{y\in\Gamma_2: x_y<y\}$, $\Gamma_2^-=\Gamma_2\setminus \Gamma_2^+$ and $\Gamma^+=\{(x,y)\in \Gamma: y\in \Gamma_2^+\}$, $\Gamma^-=\Gamma\setminus \Gamma^+$.\\
First of all note that upon stopping at $T_y$ we have that $H(y,y)=F(y)$ which means that $\Gamma$ cannot contain points from the diagonal in $\re_+^2$.\\
Let $\mathcal{R}$ be the range of the process $(B_t,\overline{B}_t)$, $\mathcal{R}(\omega)=\{(B_t(\omega),\overline{B}_t(\omega)):t\ge 0\}$.
The restriction we have to impose on $\Gamma$ is that $\p\big( \mathcal{R}\cap \Gamma=\emptyset)=1$. Notice however that, due to the continuity of sample paths of $B$, if $(x,y)\in \mathcal{R}\cap \Gamma$ then $(x_y,y)\in \mathcal{R}$. With L\'evy's equivalence theorem we know that the process $\gamma_t=(\overline{B}_t-B_t)$ is a reflected standard Brownian motion and the stretches $[x,y]\times \{y\}$ in $\mathcal{R}$ correspond to its excursions, which form a Poisson point process on the time scale given by $\overline{B}_t$. Thus the process of extremal values of these excursions, $(\overline{e}_y:y\ge 0)=(\sup_{s\le T_{y+}-T_y}(y-B_s):y\ge 0)$, is also a Poisson point process and its intensity measure is $\frac{dv}{v^2}$ (cf.\ Revuz and Yor \cite{MR2000h:60050} ex.\ XII.2.10). Then we have $\p(\mathcal{R}\cap \Gamma^+=\emptyset)=\exp(-\int_{\Gamma_2^+}\frac{dy}{y-x_y})$, which is equal to one if and only if $|\Gamma_2^+|=0$, where $|\cdot|$ denotes the Lebesgue measure on $\re$.\\
The probability $\p(\mathcal{R}\cap\Gamma^-=\emptyset)$ is just the probability that the Poisson point process $(\overline{e}_y:y\ge 0)$ has no jumps for $y\in \Gamma^-_2$ and this probability is zero if and only if $|\Gamma_2^-|=0$. Indeed if, for $A\subset \re_+^2$, we note $N(A)$ the cardinality of $\{y:(y,\overline{e}_y)\in A\}$ then we have
$\p(\mathcal{R}\cap\Gamma^-=\emptyset)=\p\big(N(\Gamma^-_2\times (0,\infty))=0\big)$ and $\p\big(N(\Gamma^-_2\times (0,\infty))>0\big)=\lim_{h\searrow 0}\p\big(N(\Gamma^-_2\times [h,\infty))>0\big)=\lim_{h\searrow 0} |\Gamma^-_2|/h$. The limit is zero if and only if $|\Gamma^-_2|=0$, which justifies our claim.\footnote{We want to thank Victor Rivero for his helpful remarks.}\\
As $\p(\mathcal{R}\cap \Gamma=\emptyset)=\p(\mathcal{R}\cap\Gamma^+=\emptyset)+\p(\mathcal{R}\cap\Gamma^-=\emptyset)$ we need to impose on $\Gamma$ that
$|\Gamma_2|=0$.
This ends the proof of Proposition \ref{prop:Hgen}.

The rest of this section is devoted to the proof of Theorem \ref{thm:maxmart} for Brownian motion.
The proof is organized in two parts. In part one we
will show that if $f:\re_+\to\re$ is a locally integrable function
and $H$ is given through (\ref{eq:Hgen}) then
$(H(B_t,\overline{B}_t):t\ge 0)$ is a local martingale, and
(\ref{eq:Hgenmartint}) holds.
In the second part we will show the converse. The first part is proved in Ob\l \'oj and Yor \cite{ja_yor_maxmart} but we quote it here for the sake of completeness.
\medskip\\
\textbf{Part 1. }
Suppose $f\in C^1$ and $H$ is given through (\ref{eq:Hgen}), so that (\ref{eq:Hgenmart}) holds by Proposition \ref{prop:Hgen}.
We can apply It\^o's formula to obtain:
\begin{eqnarray*}
  \label{eq:itoH}
  H(B_t,\overline{B}_t)&=&H(0,0)+\int_0^t f(\overline{B}_s)dB_s+\int_0^tf'(\overline{B}_s)(\overline{B}_s-B_s)d\overline{B}_s\\&=&H(0,0)+\int_0^t f(\overline{B}_s)dB_s,\quad
\textrm{since } d\overline{B}_s-\textrm{a.s.\ }B_s=\overline{B}_s.
\end{eqnarray*}
We have established thus the formula (\ref{eq:Hgenmartint}) for $f$ of class $C^1$.
Thus if we can show that the quantities given in
(\ref{eq:Hgenmart}) and (\ref{eq:Hgenmartint}) are well defined and finite for any locally
integrable $f$ on $[0,\infty)$, then the formula (\ref{eq:Hgenmart})--(\ref{eq:Hgenmartint})
extends to such functions through monotone class theorem.
In particular, we see that $(H(B_t,\overline{B}_t):t\ge 0)$, for $H$ given by (\ref{eq:Hgen}), is a local martingale, as it is a stochastic integral with respect to Brownian motion.
For $f$ a locally integrable function, $F(x)$ is well defined and finite, so all we
need to show is that $\int_0^t f(\overline{B}_s)dB_s$ is well
defined and finite a.s.\ for all $t>0$. This is equivalent to $\int_0^t
\Big(f(\overline{B}_s) \Big)^2ds<\infty$ a.s., for all $t>0$ which we now show.

Write $T_x=\inf\{t\ge 0: B_t=x\}$ for the first hitting time of $x$,
which is a well defined, a.s.\ finite, stopping time. Integrals in question are finite, $\int_0^t
\Big(f(\overline{B}_s) \Big)^2ds<\infty$ a.s., for all $t>0$, if and only if, for
all $x>0$, $\int_0^{T_x} \Big(f(\overline{B}_s) \Big)^2ds<\infty$.
However, the last integral can be rewritten as
\begin{eqnarray}
  \label{eq:localvsqlocal}
  \int_0^{T_x}ds\Big(f(\overline{B}_s)\Big)^2&=&\sum_{0\le u\le x}\int_{T_{u-}}^{T_u}ds \Big(f(\overline{B}_s)\Big)^2\nonumber\\
&=&\sum_{0\le u \le x}f^2(u)\Big(T_u-T_{u-}\Big)=\int_0^x
f^2(u)dT_u.
\end{eqnarray}
Now it suffices to note that\footnote{ Recall that $(T_x:x\ge 0)$ is a $\frac{1}{2}$-stable subordinator. Equality (\ref{eq:stabsub}) is easily established for simple functions and passage to the limit (cf.\ Revuz and Yor \cite{MR2000h:60050} pp.\ 72, 107 and Ex. III.4.5).}
\begin{equation}
\label{eq:stabsub}
  \e\Big[\exp\Big(-\frac{1}{2}\int_0^xf^2(u)dT_u\Big)\Big]=\exp\Big(-\int_0^x|f(u)|du\Big),
\end{equation}
to see that the last integral in (\ref{eq:localvsqlocal}) is finite
if and only if $\int_0^x|f(u)|du<\infty$, which is precisely our
hypothesis on $f$. Finally, note that the function $H$ given by (\ref{eq:Hgen}) is locally
integrable as both $x\to f(x)$ and $x\to xf(x)$ are locally
integrable.\smallskip\\
\textbf{Part 2. } In this part we show the converse to the first
part. Namely, we show that if $H:\mathbf{D}\to \re$ is a Borel function
such that $(H(B_t,\overline{B}_t):t\ge 0)$ is a right-continuous local martingale, then there exists a locally integrable function $f:\re_+\to\re$ such (\ref{eq:Hgenmart}) holds. (\ref{eq:Hgenmartint}) then holds by Part 1 of the proof above and $H$ is described by Proposition \ref{prop:Hgen}.
We start with a lemma.
\begin{lemm}
\label{lem:cle}
Let $r>0$ and $K:(-\infty, r]\to\re$ be a Borel function, such that
$(K(B_{t\land T_r}):t\ge 0)$ is a right-continuous local martingale.
Then there exist a constant $\alpha$ such that $K(x)=\alpha x+K(0)$ for $x\in (-\infty,r]$ and $(K(B_{t\land T_r}):t\ge 0)$ is a martingale.
\end{lemm}
\begin{proof}
This lemma essentially says that the scale functions for Brownian motion are the affine functions. This is a well known fact, however, for the sake of completeness, we provide a short proof.\\
We know that a right-continuous local martingale has actually a.s.\ c\`adl\`ag paths. Furthermore, as $K_t=K(B_{t\land T_r})$ is a local martingale with respect
to the Brownian filtration generated by $B$, it actually has a continuous version. As
the laws in the space of c\`adl\`ag functions are determined by
finite-dimensional projections, $K_t$ is a.s.\ continuous, which
implies that $K(\cdot)$ is continuous on $(-\infty,r]$.\footnote{We could also just say that the right-continuity of $K(B_t)$ implies fine-continuity of $K$ (cf.\ Thm II.4.8 in Blumenthal and Getoor \cite{MR41:9348}) and the fine topology for real-valued Brownian motion is the ordinary topology.}
\\
Let $T_{a,b}=\inf\{t\ge 0: B_t\notin [a,b]\}$. Then, as $K$ is bounded on compact sets, for any $0<x<r$ the local martingale $K(B_{t\land T_{-1,x}\land T_r})=K(B_{t\land T_{-1,x}})$ is bounded and hence it is a uniformly integrable martingale. Applying the optional stopping theorem we obtain $\e K(B_{T_{-1,x}})=K(0)$ and thus $K(x)=x(K(0)-K(-1))+K(0)$. Similarly, for $x<0$, we can apply the optional stopping theorem to see that $\e K(B_{T_{x,r/2}})=K(0)$, which yields $K(x)=x\frac{2K(r/2)-2K(0)}{r}+K(0)$. As $K$ is continuous, we conclude that it is an affine function on $(-\infty, r]$.\footnote{We thank Goran Peskir and Dmitry Kramkov for their remarks, which simplified our earlier proof of the lemma.}
\end{proof}
We now turn to the proof of the theorem. We will show how it reduces
to the above lemma. With no loss of generality we can assume that $H(0,0)=0$.
The proof is carried out in $5$ steps:
\begin{enumerate}
    \item For almost all $y$, the function $H(\cdot,y)$ is continuous.
    \item For all $y>x\lor 0$ and suitable random variables $\xi$ independent of Brownian motion $\beta$,
    $(H(x+\beta_{t\land R_{y-x}},\xi):t\ge 0)$ is a local martingale (where $R_u=\inf\{t:\beta_t=u\}$) in the filtration of $\beta$ enlarged with $\xi$.
    \item For almost all $z$, $z>y>x\lor 0$, actually $(H(x+\beta_{t\land R_{y-x}},z):t\ge 0)$ is a local martingale in the natural filtration of $\beta$.
    \item Apply Lemma \ref{lem:cle} to obtain (\ref{eq:Hgenmart}).
    \item Proof of the martingale property.
\end{enumerate}
\textbf{Step 1. }
As in the proof of Lemma \ref{lem:cle} we can argue that $(H(B_t,\overline{B}_t):t\ge 0)$ is a continuous local martingale.
From the proof of Proposition \ref{prop:Hgen} above, in particular from the discussion of the range of the process $(B,\overline{B})$, it follows that for almost all $z\ge 0$, $H(\cdot,z)$ is a continuous function on $(-\infty, z]$.
As we want to prove the a.s.\ representation given by (\ref{eq:Hgenmart}) we know, by Proposition \ref{prop:Hgen}, that we can change $H$ on some set of the form $\cup_{z\in A}(-\infty,z)\times \{z\}$ with $A$ of zero Lebesgue measure, and so
we can and will assume that $H(\cdot,z)$ is continuous on $(-\infty,z]$ for all $z\ge 0$.

Let $y>x\lor 0$ and $T_y=\inf\{t\ge 0: B_t=y\}$, and $T^y_{x}=\inf\{t>T_y:B_t=x\}$, two almost surely finite stopping times. Denote $\xi=\overline{B}_{T^y_{x}}$, which is a random variable with an absolutely continuous distribution on $[y,\infty)$. We note its density $\rho$, $\p(\xi\in du)=\rho(u)\mathbf{1}_{u\ge y}du$. We will need this notation in the sequel. Note that we could also derive the continuity properties of $H$ analyzing the behavior of $H(B_t,\overline{B}_t)$ for $t$ between the last visit to $y$ before $T^y_x$ and $T^y_x$.
\smallskip\\
\textbf{Step 2. }
Without any loss of generality we may assume $H(0,0)=0$. Using the representation theorem for Brownian martingales we know that there exists a predictable process $(h_s:s\ge 0)$ such that $H(B_t,\overline{B}_t)=\int_0^t h_s dB_s$ a.s. Let $(\theta_t:t\ge 0)$ be the standard shift operator for the two-dimensional Markov process $((B_t,\overline{B}_t):t\ge 0)$. Obviously, for $t,s>0$, we have
\begin{eqnarray}
\label{eq:addhelp2}  H(B_{t+s},\overline{B}_{t+s})-H(B_t,\overline{B}_t)&=&\Big(H(B_s,\overline{B}_s)-H(B_0,\overline{B}_0)\Big)\circ \theta_t.\nonumber
\end{eqnarray}
If we rewrite this, using the integral representation, we see that
\begin{eqnarray}
\label{eq:addhelp3}
\int_0^s h_{u+t}dB_{u+t}=\int_0^s\Big(h_u\circ \theta_t\Big)d B_{u+t}\quad a.s.,\end{eqnarray}
which implies that $h_{u+t}=h_u\circ \theta_t$ for $u>0$ a.s. Reasoning stays true if we replace $t$ by an arbitrary, a.s.\ finite, stopping time $T$. This in turn means that the process $\langle H(B,\overline{B}),B\rangle_t = \int_0^t h_s ds$ is a signed (strong) additive functional of the process $((B_t,\overline{B}_t):t\ge 0)$. To each of the strong additive functionals $\int_0^t (h_s\lor 0)ds$ and $\int_0^t (-h_s\lor 0)ds$ we can apply Motoo's theorem (cf.\ Sharpe \cite{MR958914} p.\ 309, see also
Meyer \cite{MR0231445} p.\ 122, and Ex. X.2.25 in Revuz and Yor \cite{MR2000h:60050}) to see that there exists a measurable function $h:\re\times\re_+\to
\re$ such that
$$H(B_t,\overline{B}_t)=\int_0^t h(B_s,\overline{B}_s)dB_s,\quad t\ge
0\quad a.s.$$

An application of the strong Markov property at the stopping time $T^y_x$, $y>x\lor 0$, defined in Step 1 above, yields that the process
 \begin{eqnarray}
\label{eq:martlok1}
   H\Big(x+\beta_t,\xi\lor(x+\overline{\beta}_t)\Big)-H(x,\xi)&=&H(B_{T^y_x+t},\overline{B}_{T^y_x+t}) -H(B_{T^y_x},\overline{B}_{T^y_x})\nonumber \\
   &=&
   \int_{0}^{t}h\Big(x+\beta_s,\xi\lor(x+\overline{\beta}_s)\Big)d\beta_s
 \end{eqnarray}
is a local martingale in the enlarged filtration
$\mathcal{G}_t=\sigma(\xi,\,\beta_s:s\leq t)$, where
$\beta_s=B_{T^y_x+s}-B_{T^y_x}$ is a new Brownian motion
independent of $(B_u:u\leq T^y_x)$. Furthermore, as on the
interval $[0,R_{y-x}]$, where $R_{y-x}=\inf\{t\ge 0:\beta_t=y-x\}$,
we have $\xi\lor (x+\overline{\beta}_s)=\xi$, the stopped local martingale can be written as
\begin{eqnarray}
\label{eq:postac_zxi}
H(x+\beta_{t\land R_{y-x}},\xi)&=&H(x,\xi)+\int_0^{t\land
R_{y-x}}h(x+\beta_t,\xi)d\beta_s.
\end{eqnarray}\smallskip\\
\textbf{Step 3. } We want to show that actually, for almost all $z\in (y,\infty)$,  $(H(x+\beta_{t\land R_{y-x}},z):t\ge 0)$ is a local
martingale in the natural filtration of $\beta$.

Let $\tilde{H}(x,z)=\big(H(x,z)-H(0,z)\big)\mathbf{1}_{x\le z}+\big(H(z,z)-H(0,z)\big)\mathbf{1}_{x>z}$, which is a measurable function, continuous in the first coordinate. Fix $K>0$ and define the function $\epsilon :(0,\infty)\to [0,1]$ via
\begin{eqnarray}
\label{eq:defep}
\epsilon (z) &=& \sup\Big\{0\leq\delta\le 1: |\tilde{H}(x,z)|\leq K \textrm{ for }x\in [-\delta K,\delta K]\Big\}\\
&=& \sup\Big\{0\leq\delta\le 1: |\tilde{H}(x,z)|\leq K \textrm{ for }x\in [-\delta K,\delta K]\cap \mathbb{Q}\Big\},
\end{eqnarray}
where the equality follows from continuity properties of $\tilde{H}$.
We now show that $\epsilon(\cdot)$ is a measurable function. To this end let $\delta\in (0,1]$ and write \begin{eqnarray*}
  \{z:\epsilon(z)<\delta\}  &=&\bigg\{z: \sup\Big\{|\tilde{H}(x,z)|:x\in [-\delta K,\delta K]\cap \mathbb{Q}\Big\}>K\bigg\} \\
 &=&\bigcup_{x\in [-\delta K,\delta K]\cap \mathbb{Q}}\Big\{z:|\tilde{H}(x,z)|>K\Big\}.
 \end{eqnarray*}
Measurability of $\epsilon$ follows as $$\Big\{z:|\tilde{H}(x,z)|>K\Big\}=\Big(\{x\}\times\re_+\Big)\cap \tilde{H}^{-1}\Big[(-\infty,K)\cup (K,\infty)\Big],$$ is a Borel set.\\
We defined $\epsilon$ so that $|\tilde{H}(x,z)|\leq K$ on
$[-\epsilon(z)K,\epsilon(z)K]$ (and it is the biggest such interval).
Note that, since a continuous function is bounded on compact intervals, we have $\epsilon(z) K\to\infty$ as $K\to\infty$.
Let $T_K$ be a stopping time in the enlarged
filtration $(\mathcal{G}_t)$ defined by $T_K=T_K(\beta,\xi)=\inf\{t\ge
0: |\beta_t|\ge \epsilon(\xi)K\}$, and write $T_K^z$ for $T_K(\beta,z)$,
which is a stopping time in the natural filtration of $\beta$.
Then
by (\ref{eq:postac_zxi}) we see that $(H(x+\beta_{t\land R_{y-x}\land
T_K},\xi)-H(x,\xi):t\ge 0)$ is a a.s.-bounded local martingale, hence a
martingale. Recall that $\rho$ is the density function of the
distribution of $\xi$, $\p(\xi\in dz)=\rho(z)\mathbf{1}_{z\ge y}dz$,
and let $b:\re\to\re$ be a Borel, bounded function. Put $x=0$. Then
the process $(M^b_t=b(\xi)H(\beta_{t\land R_y\land T_K},\xi):t\ge 0)$
is a $(\mathcal{G}_t)$-martingale. We want to show that the process $H(\beta_{t\land T_K^z\land R_y},z)$ is a martingale for almost all $z>y$. As we deal with continuous, a.s.-bounded processes in a continuous filtration, it suffices to verify the martingale property for rational times.
Fix $t,s\in \mathbb{Q}_+$. For any $A\in
\mathcal{F}^\beta_t=\sigma(\beta_s:s\leq t)$, by the martingale property of $M^b$, we have
\begin{eqnarray}
&&  \e\Big[\mathbf{1}_A b(\xi)H(\beta_{(t+s)\land R_y\land T_K},\xi)\Big]=   \e\Big[\mathbf{1}_A b(\xi)H(\beta_{t\land R_y\land T_K},\xi)\Big]\nonumber \\
 && \int_y^\infty dz\rho(z)b(z) \e\Big[\mathbf{1}_A H(\beta_{(t+s)\land R_y\land T^z_K},z)\Big]\nonumber\\
 &&\qquad\qquad=\int_y^\infty dz\rho(z)b(z) \e\Big[\mathbf{1}_A H(\beta_{t\land R_y\land T^z_K},z)\Big] \textrm{ and as }b\textrm{ was arbitrary},\nonumber\\
  &&\e\Big[\mathbf{1}_A H(\beta_{(t+s)\land R_y\land T^z_K},z)\Big]= \e\Big[\mathbf{1}_A H(\beta_{t\land R_y\land
  T^z_K},z)\Big],\quad dz-a.e.\label{eq:marthp}
\end{eqnarray}
We will now argue that the above actually holds $dz$-a.e.\ for all $t,s\in\mathbb{Q}$ and all $A\in \mathcal{F}^\beta_t$.
Let $\Pi\subset \mathcal{F}^\beta_t$ be a countable $\pi$-system which generates $\mathcal{F}_t^\beta$ (cf.\ Exercise I.4.21 in Revuz and Yor \cite{MR2000h:60050}). We can thus choose a set $\Gamma_{t,s}\subset (y,\infty)$ of full Lebesgue measure, such that for any $A\in \Pi$, (\ref{eq:marthp}) holds for all $z\in \Gamma_{t,s}$. As sets $A$ which satisfy (\ref{eq:marthp}) for all $z\in \Gamma_{t,s}$ form a $\lambda$-system, it follows that (\ref{eq:marthp}) holds for any $A\in \mathcal{F}^\beta_t$ for all $z\in \Gamma_{t,s}$. Letting $\Gamma=\bigcap_{t,s\in \mathbb{Q}_+}\Gamma_{t,s}$, we see that $(H(\beta_{t\land T_K^z\land R_y},z):t\ge 0)$ is a martingale for all $z\in \Gamma$, and $\Gamma$ is of full Lebesgue measure. This implies the local martingale property for
$(H(\beta_{t\land R_y},z):t\ge 0)$ since $(T_K^z:K\in \nr)$ is a good localizing sequence. Indeed, for almost all $z>y$, $\epsilon(z)K\to \infty$ as $K\to \infty$, and so
 $T_K^z\to\infty$ a.s., as $K\to \infty$.\smallskip\\
\textbf{Step 4. } We know thus that for almost all $z>y$,
$(H(\beta_{t\land R_{y}},z):t\ge 0)$ is a local martingale with respect to
the natural filtration of $\beta$. We can thus apply Lemma \ref{lem:cle}
 to see that $K(b)=H(b,z)$ is a linear function on $(-\infty, y]$, for almost all $z>y$. Thus $H(\beta_{t\land
R_{y}},\xi)=\alpha (\xi) \beta_{t\land R_{y}}+H(0,\xi)$ a.s. Confronting this
with (\ref{eq:postac_zxi}) we see that $h(b,z)$ does not depend on
$b$ for $b\in (-\infty,y    ]$, $h(b,z)=h(z)$ for almost all $z>y$.
As $y>0$ was arbitrary, taking $y\in \mathbb{Q}$ and $y\to 0$, we see that $h(b,z)=h(z)$ for almost all
$z>0$, and therefore $H(B_t,\overline{B_t})=\int_0^t
h(\overline{B}_u)dB_u$ a.s. and we put $f=h$.
From the first part of the proof we know that if $\int_0^t f(\overline{B}_s)dB_s$ is
well defined and finite then $f$ is locally integrable and (\ref{eq:Hgenmart}) holds.
\smallskip\\
\textbf{Step 5. }
We turn now to the proof of the last statement in Theorem \ref{thm:maxmart}. Let $f$ be a Borel, positive function in $L^1$, $||f||=\int_0^\infty f(x)dx$.
Define $H(x,y)=||f||-F(y)+f(y)(y-x)$. The process $H(N_t,\overline{N}_t)$ is a local martingale as in (\ref{eq:Hgenmart}). Furthermore, it is a positive process and we can apply Fatou's lemma to see that it is a positive supermartingale and thus converges a.s., as $t\to\infty$. As $\langle N\rangle_\infty=\infty$ we know that $\overline{N}_\infty=\infty$ a.s. This entails that $||f||-F(\overline{N}_t)=\int_{\overline{N}_t}^\infty f(x)dx\vat 0$ a.s.\ and thus $f(\overline{N}_t)(\overline{N}-N_t)$ also converges a.s.\ as $t\to\infty$. However it can only converge to zero since $\overline{N}_\infty=\infty$ a.s.\ and thus $(\overline{N}_t-N_t)$ has zeros for arbitrary large $t$. Convergences announced in Theorem \ref{thm:maxmart} and in the remarks which followed it are immediate.\\
To establish the martingale property it suffices to see that the expectation of the positive supermartingale $H(N_t,\overline{N}_t)$ is constant in time and equal to $H(0,0)=||f||$. If $f$ is bounded and $\e(\sup_{s\le t} |N_s|)<\infty$, this follows readily from Lebesgue's dominated convergence theorem. The general case follows from monotone convergence theorem by replacing $f$ with $\min\{f,n\}$ and taking the limit as $n\to\infty$. This ends the proof of Theorem \ref{thm:maxmart}.
\section{Optional stopping arguments}
\label{sec:optional}
In the previous section we proved Theorem \ref{thm:maxmart}. Here we want to present some alternative arguments which could have been used in the proof and which relay on the optional stopping theorem.\\
We place ourselves in Brownian motion setup, that is $N=B$ is a real-valued Brownian motion.
Note that the characterization of max-martingales proved in Theorem \ref{thm:maxmart} justifies the application of the optional stopping theorem  to the local martingale displayed in (\ref{eq:Hgenmart}) at the first exit time of the underlying Brownian motion from an interval. If one could justify this independently, the following reasoning could replace some parts of our proof of Theorem \ref{thm:maxmart}.\\
As above, we may assume $H(0,0)=0$.
Let $T_{x,y}=\inf\{t\ge 0: B_t\notin [x,y]\}$ and fix $x,a,y$ in $\re_+$. Recall that we argued that $H(\cdot, y)$ is continuous on $(-\infty,y]$.
Relaying on the first part of the proof of Theorem \ref{thm:maxmart}, which grants that
the processes defined via (\ref{eq:Hgenmart}) are local martingales, one can verify the well known fact that the law of $\overline{B}_{T_{-x,y}}\mathbf{1}_{\{B_{T_{-x,y}}=-x\}}$ is given by
$\p(\overline{B}_{T_{-x,y}}\mathbf{1}_{\{B_{T_{-x,y}}=-x\}}\in ds)=\frac{xds}{(s+x)^2}\mathbf{1}_{0\le s\le y}$. An application of the optional stopping theorem to $H(B_t,\overline{B}_t)$ at $T_{-x,y}$ and $T_{-x-a,y}$ yields:
\begin{eqnarray}
  0 &=& \frac{x}{x+y}H(y,y)+ x\int_0^y\frac{H(-x,s)}{(x+s)^2}ds,\quad \textrm{and}\\
0 &=& \frac{x+a}{x+a+y}H(y,y)+ (x+a)\int_0^y\frac{H(-x-a,s)}{(x+a+s)^2}ds.
\end{eqnarray}
Solving both equations for $H(y,y)$ and comparing leads to
\begin{equation}\label{eq:1}
    (x+y)\int_0^y\frac{H(-x,s)}{(x+s)^2}ds=(x+a+y)\int_0^y\frac{H(-x-a,s)}{(x+a+s)^2}ds.
\end{equation}
Both sides are differentiable in $y$ and differentiating we obtain
\begin{eqnarray}
\label{eq:2}
\int_0^y\frac{H(-x,s)}{(x+s)^2}ds+\frac{H(-x,y)}{x+y}   &=&  \int_0^y\frac{H(-x-a,s)}{(x+a+s)^2}ds+\frac{H(-x-a,y)}{x+a+y}\quad\textrm{ thus}\nonumber\\
 \frac{H(-x-a,y)}{x+a+y}-\frac{H(-x,y)}{x+y}&=&\int_0^y\frac{H(-x,s)}{(x+s)^2}ds-\int_0^y\frac{H(-x-a,s)}{(x+a+s)^2}ds\nonumber\\
&=&\frac{a}{x+a+y}\int_0^y\frac{H(-x,s)}{(x+s)^2}ds,\quad \textrm{by }(\ref{eq:1}).
\end{eqnarray}
Transforming the last equality we obtain finally
\begin{equation}\label{eq:3}    H(-x-a,y)=a\Big(\int_0^y\frac{H(-x,s)}{(x+s)^2}ds+\frac{H(-x,y)}{x+y}\Big)+H(-x,y),
\end{equation}
and letting $x\to 0$, as $H(\cdot, y)$ is continuous, we have
\begin{equation}\label{eq:3}
    H(-a,y)=a\Big(\int_0^y\frac{H(0,s)}{s^2}ds+\frac{H(0,y)}{y}\Big)+H(0,y),
\end{equation}
which shows that $H(\cdot,y)$ is a linear function on $(-\infty,0)$. Furthermore, if we define $F(y)=y\int_0^y\frac{H(0,s)}{s^2}ds$ then $F(0)=0$ and $f(y)=F'(y)=\int_0^y\frac{H(0,s)}{s^2}ds+\frac{H(0,y)}{y}$. We can now rewrite (\ref{eq:3}) as
\begin{equation}\label{eq:4}
    H(-a,y)=(a+y)f(y)-F(y),\quad a,y>0,
\end{equation}
in which we instantly recognize the desired form displayed in (\ref{eq:Hgen}).\\
However, in order to recover the desired form of the function $H$ on the whole set $\mathbf{D}$ we would need to argue that for any $u>0$, $(H(u+B_t,u+\overline{B}_t):t\ge 0)$ is also a local martingale. This follows from the Step 2 in our proof of the second part of Theorem \ref{thm:maxmart} but might not be easy to see independently.

\section{Closing remarks}
At first glance Theorem \ref{thm:maxmart} has mainly the theoretical value of providing a complete characterization of a certain family of local martingales. To close this paper, we point out that it has some further interesting consequences.

A useful method of proving various inequalities, such as Doob-like inequalities, consists in exhibiting appropriate martingales and applying the optional stopping theorem
(cf.\ Ob\l \'oj and Yor \cite{ja_yor_maxmart}). Theorem \ref{thm:maxmart} tells us that if we search for a martingale which involves only Brownian motion and its maximum process then we have to look among the Az\'ema-Yor martingales.
As mentioned above, this work is a first step in a more general project of describing local martingales which are function of Brownian motion and some adapted, $\re^d$-valued, process with ``small support'' (as the maximum, minimum and local time processes). Such martingales, for functions which are regular enough, can be described via It\^o's formula. However, a complete characterization for arbitrary functions is more delicate. We believe that the methodology developed in our proof of Theorem \ref{thm:maxmart} will be useful for this purpose.

The second remark we want to make is in close link with some penalization problems discussed by Roynette, Vallois and Yor \cite{rvy2}. They remark in \cite{rvy2} that the limiting martingales they obtain have a special form. It seems that this can be justified with a similar application of Motoo's Theorem as in the proof of Theorem \ref{thm:maxmart} above.
We plan to develop this topic in a separate paper.\smallskip\\
\textbf{Acknowledgment.} I am deeply indebted to Marc Yor whose ideas and help were essential for the development of this paper.
\bibliography{/users/obloj/docs/Moje_dokumenty/Matematyka/bib/bibliografia.bib}
\end{document}